# Sets of Completely Decomposed Primes in Extensions of Number Fields

by Kay Wingberg

Let $p$ be a prime number and let $k(p)$ be the maximal $p$-extension of a number field $k$. If $T$ is a set of primes of $k$, then $k^T(p)$ denotes the maximal $p$-extension of $k$ which is completely decomposed at $T$. Assuming that $T$ is finite, the canonical homomorphism

$$\phi^T(p): \underset{\mathfrak{p}\in T(k^T(p))}{\bigstar} G_{\mathfrak{P}}(k(p)|k) \longrightarrow G(k(p)|k^T(p))$$

of the free pro-$p$-product of the decomposition groups $G_{\mathfrak{P}}(k(p)|k)$ into the Galois group $G(k(p)|k^T(p))$ is an isomorphism, see [3] theorem (10.5.8); here the prime $\mathfrak{P}$ is an arbitrary extension of $\mathfrak{p}$ to $k(p)$.

In the profinite case, i.e. considering the maximal Galois extension $k^T$ which is completely decomposed at the finite set $T$, there exist suitable extensions $\mathfrak{P}|\mathfrak{p}$ such that $\phi^T$ is an isomorphism of profinite groups. If $T = S_\infty$ is the set archimedean primes, this is a result of Fried-Haran-Völklein [1] and in general it is proven by Pop [4].

The fact that $\phi^T(p)$ is an isomorphism if $T$ is finite implies very strong properties for the extension $k^T(p)|k$. In particular, the (strict) cohomological dimension of the Galois group $G(k^T(p)|k)$ is equal to 2 (if $p = 2$ one has to require that $k$ is totally imaginary). Furthermore, we get for the corresponding local extensions that

$(*)$  $(k^T(p))_{\mathfrak{P}} = k_{\mathfrak{p}}(p)$ for all primes $\mathfrak{P}|\mathfrak{p}$, $\mathfrak{p} \notin T$,

i.e. $k^T(p)$ realizes the maximal $p$-extension $k_{\mathfrak{p}}(p)$ of the local fields $k_{\mathfrak{p}}$ for all primes not in $T$. In particular, the set $T$ is equal to the set $D(k^T(p)|k)$ of all primes of $k$ which decomposed completely in the extension $k^T(p)|k$. We will say that $T$ is *saturated* if it has this property and we call the stronger property $(*)$ that $T$ is *strongly saturated*. If the Dirichlet density $\delta(T)$ is positive, then $T$ is saturated if and only if $T$ is the set of completely decomposed primes of a *finite* Galois extension of $k$.

If $T$ is an arbitrary set of primes of $k$, we call the set $\hat{T} = D(k^T(p)|k)$ the *saturation* of $T$. The most interesting case is that $T$ is an infinite set of primes of density zero. If $T$ is saturated, then the extension $k^T(p)|k$ is infinite. We will see that there exist infinite sets $T$ of primes such that $\delta(\hat{T}) > \delta(T) = 0$, and also sets $T$ such that $\delta(\hat{T}) = \delta(T) = 0$. An important example is the following.

**Theorem 1:** *Let $p$ be an odd prime number and let $k$ be a CM-field containing the group $\mu_p$ of all $p$-th roots of unity, with maximal totally real subfield $k^+$, i.e. $k = k^+(\mu_p)$ is totally imaginary and $[k : k^+] = 2$. Let $S_p = \{\mathfrak{p}|p\}$ and*

$$T = \{\mathfrak{p} \mid \mathfrak{p} \cap k^+ \text{ is inert in } k|k^+\}.$$

*Then $T \cup S_p$ is strongly saturated.*

Moreover, in the example above we get that the Galois group

$$G((k^T)_{nr}(p)|k^T(p))$$

of the maximal unramified $p$-extension of $(k^T)_{nr}(p)$ of $k^T(p)$ is a free pro-$p$-group. This will follow from a more general theorem which deals with a generalization of the notion of saturated sets.

A set $T = T(k)$ is called *stably saturated* if the sets $T(k')$ are saturated for every finite Galois extension $k'|k$ inside $k^T(p)$. These sets are necessarily of density 0 if they are not equal to set $\mathcal{P}$ of all primes. Obviously, strongly staturated sets are stably saturated. We have the following theorem.

**Theorem 2:** *Let $T \neq \mathcal{P}$ be a stably saturated set of primes of a number field $k$. Then the canonical map*

$$\phi^T(p) : \underset{\mathfrak{P} \in T(k^T(p))}{\bigstar} (G_\mathfrak{P}(k(p)|k), T_\mathfrak{P}(k(p)|k)) \xrightarrow{\sim} G(k(p)/k^T(p))$$

*is an isomorphism.*

Here $\bigstar_{\mathfrak{P} \in T(k^T(p))}(G_\mathfrak{P}(k(p)|k), T_\mathfrak{P}(k(p)|k))$ denotes the free corestricted pro-$p$-product of the decomposition groups $G_\mathfrak{P}(k(p)|k)$ with respect to the inertia groups $T_\mathfrak{P}(k(p)|k)$, see [2].

**Corollary:** *In the situation of theorem 1 let $\tilde{T} = T \cup S_p$. Then*

$$\underset{\mathfrak{P} \in \tilde{T}(k^{\tilde{T}}(p))}{\bigstar} (G_\mathfrak{P}(k(p)|k), T_\mathfrak{P}(k(p)|k)) \xrightarrow{\sim} G(k(p)/k^{\tilde{T}}(p)).$$

Finally I would like to thank Jochen Gärtner for helpful discussions on this subject and Alexander Ivanov for valuable comments.

# 1  Saturated sets of primes of a number field

We start with some remarks on complete lattices. Let

$$\mathcal{A} \underset{\psi}{\overset{\varphi}{\rightleftarrows}} \mathcal{B}$$

be maps between complete lattices $(\mathcal{A}, \subseteq)$ and $(\mathcal{B}, \subseteq)$ with the following properties:

I. $\varphi$ and $\psi$ are order-reversing,

II. $A \subseteq \psi\varphi(A)$ and $B \subseteq \varphi\psi(B)$ for all $A \in \mathcal{A}$ and $B \in \mathcal{B}$.

For $A \in \mathcal{A}$ and $B \in \mathcal{B}$ we define the *saturation*

$$\hat{A} := \psi\varphi(A) \text{ and } \hat{B} := \varphi\psi(B),$$

and call $A \in \mathcal{A}$ resp. $B \in \mathcal{B}$ to be *saturated* if $A = \hat{A}$ resp. $B = \hat{B}$. We put

$$\mathcal{A}_{sat} = \{A \in \mathcal{A} \mid A \text{ is saturated}\}, \qquad \mathcal{B}_{sat} = \{B \in \mathcal{B} \mid B \text{ is saturated}\},$$

and we have the following properties:

(i) $A_1 \subseteq A_2$ implies $\hat{A}_1 \subseteq \hat{A}_2$ and $B_1 \subseteq B_2$ implies $\hat{B}_1 \subseteq \hat{B}_2$.

(ii) $\varphi\psi\varphi = \varphi$, $\psi\varphi\psi = \psi$, $\hat{\hat{A}} = \hat{A}$, $\hat{\hat{B}} = \hat{B}$.

(iii) $\mathcal{B}_{sat}$ is the image of $\mathcal{A}$ under $\varphi$ and $\mathcal{A}_{sat}$ is the image of $\mathcal{B}$ under $\psi$, i.e. $\varphi: \mathcal{A} \twoheadrightarrow \mathcal{B}_{sat}$ and $\psi: \mathcal{B} \twoheadrightarrow \mathcal{A}_{sat}$, and $\psi$ and $\varphi$ induce bijections

$$\mathcal{A}_{sat} \underset{\psi}{\overset{\varphi}{\rightleftarrows}} \mathcal{B}_{sat}.$$

(iv)
$$\varphi(\bigcup_i A_i) = \bigcap_i \varphi(A_i), \qquad \varphi(\bigcap_i \psi(B_i)) = \widehat{\bigcup_i B_i}.$$

In particular, the infimum of saturated elements is again saturated. The verification of these statements is straightforward using the properties I and II of the maps $\varphi$ and $\psi$. We define the equivalence relations

$$A_1 \sim A_2 \quad :\Leftrightarrow \quad \varphi(A_1) = \varphi(A_2) \qquad \text{for} \quad A_1, A_2 \in \mathcal{A}$$

and

$$B_1 \sim B_2 \quad :\Leftrightarrow \quad \psi(B_1) = \psi(B_2) \qquad \text{for} \quad B_1, B_2 \in \mathcal{B},$$

and denote the classes by $[X]$. Obviously, $X \sim \hat{X}$ by (ii) and for $Y \in [X]$ we have $Y \subseteq \hat{X}$, and so $\hat{X}$ is the unique maximal element of $[X]$. Furthermore, the surjections $\psi\varphi: \mathcal{A} \twoheadrightarrow \mathcal{A}_{sat}$ and $\varphi\psi: \mathcal{B} \twoheadrightarrow \mathcal{B}_{sat}$ induce bijections

$$\psi\varphi: (\mathcal{A}/_\sim) \xrightarrow{\sim} \mathcal{A}_{sat}, \quad [A] \mapsto \hat{A}, \quad \text{and} \quad \varphi\psi: (\mathcal{B}/_\sim) \xrightarrow{\sim} \mathcal{B}_{sat}, \quad [B] \mapsto \hat{B}.$$

Now let $K$ be a number field. We use the following notation: $S_\infty$, $S_\mathbb{R}$ and $S_\mathbb{C}$ are the sets of archimedean, real and complex primes of $K$, respectively, $\mathcal{P}$ is the set of all primes of $K$, and if $p$ is a prime number, then $S_p$ is the set of all primes of $K$ above $p$.

If $\mathfrak{p}$ is a prime of $K$, then $K_\mathfrak{p}$ is the completion of $K$ with respect to $\mathfrak{p}$. If $L|K$ is a Galois extension, then we denote the decomposition group and inertia group of the Galois group $G(L|K)$ with respect to $\mathfrak{p}$ by $G_\mathfrak{p}(L|K)$ and $T_\mathfrak{p}(L|K)$, respectively.

For a set $S$ of primes of $K$, let $\delta(S) = \delta_K(S)$ be its Dirichlet density. If $S = S(K)$ is a set of primes and $K'|K$ an algebraic extension of $K$, then we denote the set of primes of $K'$ consisting of all prolongations of $S$ by $S(K')$.

If $\mathfrak{c}$ is a class of finite groups which is closed under taking subgroups, homomorphic images and group extensions, then $K(\mathfrak{c})$ is the maximal pro-$\mathfrak{c}$-extension of $K$, and in particular if $p$ is a prime number, $K(p)$ denotes the maximal $p$-extension of $K$. By abuse of notation, we denote the maximal pro-$\mathfrak{c}$ extension of $K$ which is completely decomposed at $T$ by $K^T(\mathfrak{c})$, and $K^T(p)$ is the maximal $p$-extension of $K$ inside $K^T$.

Let
$$\mathcal{E}_K = \{L \mid L \text{ is a Galois extension of } K\} \xrightleftharpoons[\psi]{\varphi} \{T \mid T \text{ is a set of primes of } K\} = \mathcal{S}_K$$

where

$\varphi(L) = D(L|K)$ is the set of primes which are completely decomposed in $L|K$,
$\psi(T) = K^T$      is the maximal Galois extension of $K$ which is completely decomposed at $T$.

Obviously, the maps $\varphi$ and $\psi$ are order-reversing and
$$L \subseteq \psi\varphi(L) = K^{D(L|K)} =: \hat{L} \quad \text{and} \quad T \subseteq \varphi\psi(T) = D(K^T|K) =: \hat{T}.$$

Furthermore, we define the equivalence relations
$$L_1 \sim L_2 \quad :\Leftrightarrow \quad D(L_1|K) = D(L_2|K) \qquad \text{for } L_1, L_2 \in \mathcal{E}_K$$

and
$$T_1 \sim T_2 \quad :\Leftrightarrow \quad K^{T_1} = K^{T_2} \qquad \text{for } T_1, T_2 \in \mathcal{S}_K.$$

**Definition 1.1** *The extension $L \in \mathcal{E}_K$ is called* **saturated** *if $\hat{L} = L$, i.e.*
$$L = K^{D(L|K)},$$

*and the set $T \in \mathcal{S}_K$ is called* **saturated** *if $\hat{T} = T$, i.e.*
$$T = D(K^T|K).$$

We strengthen the notion of saturated sets in the following way:

**Definition 1.2** *Let $T$ be a set of primes of $K$.*

(i) *The set $T = T(K)$ is called* **stably saturated** *if the sets $T(K') \in \mathcal{S}_{K'}$ are saturated for every finite Galois extension $K'|K$ inside $K^T$.*

(ii) *We call $T$ to be* **strongly saturated** *if*

$$(K^T)_{\mathfrak{P}} = \bar{K}_{\mathfrak{p}} \quad \text{for all primes} \quad \mathfrak{P}|\mathfrak{p},\ \mathfrak{p} \notin T,$$

*where $\bar{K}_{\mathfrak{p}}$ is the algebraic closure of $K_{\mathfrak{p}}$.*

**Remark 1:** If we consider the set $\mathcal{E}_K(\mathfrak{c}) = \{L \mid L \text{ is a pro-}\mathfrak{c}\text{-extension of } K\}$, we define a set $T \in \mathcal{S}_K$ to be $\mathfrak{c}$-saturated if $T = D(K^T(\mathfrak{c})|K)$, and an extension $L \in \mathcal{E}_K(\mathfrak{c})$ is called $\mathfrak{c}$-**saturated** if $L = K^{D(L|K)}(\mathfrak{c})$. Analogously, we define stably $\mathfrak{c}$-saturated and strongly $\mathfrak{c}$-saturated, e.g. $T$ is strongly $\mathfrak{c}$-saturated if $(K^T(\mathfrak{c}))_{\mathfrak{P}} = K_{\mathfrak{p}}(\mathfrak{c})$ for all primes $\mathfrak{P}|\mathfrak{p}$, $\mathfrak{p} \notin T$.

We say, a prime of $K$ is **redundant** (or more precise, $\mathfrak{c}$-redundant) if is totally decomposed in every extension inside $K(\mathfrak{c})$. Obviously, redundant primes are necessarily archimedean primes and the complex primes are always redundant. But also real primes might be redundant if we restrict to pro-$\mathfrak{c}$-extensions, e.g. if we consider $p$-extensions where $p$ is an odd prime number. Therefore we make the following

**Convention:** In the following all considered primes are not redundant and $\mathcal{S}_K = \mathcal{S}_K(\mathfrak{c})$ consists only of sets of non-redundant primes for the extension $K(\mathfrak{c})|K$.

Immediately from the definitions (and the convention) above we get

**Lemma 1.3** *Let $T \in \mathcal{S}_K$.*

(i) *$T$ is saturated if and only if $(K^T)_{\mathfrak{p}} \neq K_{\mathfrak{p}}$ for all $\mathfrak{p} \notin T$.*

(ii) *$T$ is stably saturated if and only if $(K^T)_{\mathfrak{p}}|K_{\mathfrak{p}}$ is an infinite extension for all $\mathfrak{p} \notin T$.*

From the general remarks on partial ordered sets we see that there are bijections

$$(\mathcal{E}_K)_{sat} \underset{\psi}{\overset{\varphi}{\rightleftarrows}} (\mathcal{S}_K)_{sat}$$

and we have the following

**Lemma 1.4**

(i) $T_1 \subseteq T_2$ implies $\hat{T}_1 \subseteq \hat{T}_2$ and $L_1 \subseteq L_2$ implies $\hat{L}_1 \subseteq \hat{L}_2$.

(ii) $K^T = K^{\hat{T}}$ and $D(L|K) = D(\hat{L}|K)$.

(iii) $D(\prod_i L_i|K) = \bigcap_i D(L_i|K)$ and $K^{\bigcup_i T_i} = \bigcap_i K^{T_i}$.

(iv) $D(\bigcap_i K^{T_i}|K) = \widehat{\bigcup_i T_i}$ and $K^{\bigcap_i D(L_i|K)} = \widehat{\prod_i L_i}$.

**Theorem 1.5**

(i) *If $T$ is a finite set of primes of $K$, then $T$ is strongly saturated.*

(ii) *If $L$ is a finite Galois extension of $K$, then $L$ is saturated.*

**Proof:** (i) Let $p$ be a prime number and let $L|K$ be a finite Galois extension inside $K^T$. Let $\mathfrak{p}_0 \notin T$, $\mathfrak{P}_0$ a fixed extension of $\mathfrak{p}_0$ to $K^T$ and $\overline{\mathfrak{P}}_0$ the restriction of $\mathfrak{P}_0$ to $L$. By the theorem of Grunwald/Wang (see [3], theorem (9.2.2)) the canonical homomorphism

$$H^1(L, \mathbb{Z}/p\mathbb{Z}) \longrightarrow H^1(L_{\overline{\mathfrak{P}}_0}, \mathbb{Z}/p\mathbb{Z}) \oplus \bigoplus_{\mathfrak{P} \in T(L)} H^1(L_{\mathfrak{P}}, \mathbb{Z}/p\mathbb{Z})$$

is surjective. In particular, for every $\alpha_{\overline{\mathfrak{P}}_0} \in H^1(L_{\overline{\mathfrak{P}}_0}, \mathbb{Z}/p\mathbb{Z}))$ there exists an element $\beta \in H^1(L, \mathbb{Z}/p\mathbb{Z})$ which is mapped to $(\alpha_{\overline{\mathfrak{P}}_0}, 0, \ldots, 0)$. But $\beta$ lies in the subgroup $H^1(K^T|L, \mathbb{Z}/p\mathbb{Z})$ of $H^1(L, \mathbb{Z}/p\mathbb{Z})$. Therefore

$$H^1(K^T|L, \mathbb{Z}/p\mathbb{Z}) \longrightarrow H^1(L_{\overline{\mathfrak{P}}_0}, \mathbb{Z}/p\mathbb{Z})$$

is surjective. Varying $L$ and $p$, it follows that the completion of $K^T$ with respect to the prime $\mathfrak{P}_0$, $\mathfrak{P}_0|\mathfrak{p}_0$ and $\mathfrak{p}_0 \notin T$, is equal to the algebraic closure of $K_{\mathfrak{p}_0}$ (since $G(\overline{K}_{\mathfrak{p}_0}|K_{\mathfrak{p}_0})$ is pro-solvable).

(ii) Let $L'$ be a finite Galois extension of $K$ with $L \subseteq L' \subseteq K^{D(L|K)}$. Since $D(L|K) \subseteq D(L'|K)$, we obtain for the densities of these sets the inequality $\delta(D(L|K)) \leq \delta(D(L'|K))$, and so, by Čebotarev's density theorem,

$$[L' : K] = \delta(D(L'|K))^{-1} \leq \delta(D(L|K))^{-1} = [L : K].$$

This shows that $L' = L$ and so $L = K^{D(L|K)}$. □

Now we consider sets of primes which are infinite and of density equal to 0.

**Proposition 1.6** *For a set $T$ of primes we have*

$$\delta(\hat{T}) = 0 \quad \Leftrightarrow \quad K^T|K \text{ is an infinite extension.}$$

**Proof:** Since $\hat{T} = D(K^T|K)$, it follows that $\delta(\hat{T}) > 0$ if $K^T|K$ is a finite extension. Conversely, assume that $K^T|K$ is infinite and let $L$ be a finite Galois extension of $K$ inside $K^T$. Then

$$\frac{1}{[L:K]} = \delta(D(L|K)) \geq \delta(\hat{T}),$$

and so $\delta(\hat{T}) = 0$. □

**Remark 2:** From the proposition above and lemma (1.3)(ii) it follows that a stably saturated set of primes $T \neq \mathcal{P}$ has necessarily density equal to 0.

**Remark 3:** Let $n \in \mathbb{N}$. Then there exist infinite sets $T$ of primes such that

$$\frac{1}{n} = \delta(\hat{T}) > \delta(T) = 0.$$

*Example 1:* Let $L|K$ be a Galois extension of the number field $K$ such that $[L:K] = n$. The set $\mathcal{E}^{fin}_{L|K}$ of proper finite extensions of $L$ being Galois over $K$ is countable, say $\mathcal{E}^{fin}_{L|K} = \{L_i, i \in \mathbb{N}\}$. We choose for every $i \in \mathbb{N}$ an element $\sigma_i \in G(L_i|L)$, $\sigma_i \neq 1$, and a prime $\mathfrak{p}_i$ of $K$ which is unramified in $L_i|K$ having a Frobenius $\left(\frac{L_i|K}{\mathfrak{P}_i}\right) = \sigma_i$, $\mathfrak{P}_i|\mathfrak{p}_i$. Since there exist infinitely many such primes, there is a $\mathfrak{p}_i$ such that $N_{K|\mathbb{Q}}\mathfrak{p}_i \geq i^2$. Then

$$\sum_{i=1}^{\infty} \frac{1}{N_{K|\mathbb{Q}}\mathfrak{p}_i} \leq \sum_{i=1}^{\infty} \frac{1}{i^2}$$

converges, and so the set $T = \{\mathfrak{p}_i, i \in \mathbb{N}\}$ has density equal to 0. Furthermore we have

$$K^T = L.$$

Indeed, since every $\mathfrak{p} \in T$ is completely decomposed in $L$, we have $L \subseteq K^T$ and a finite Galois extension $E|K$, $L \subsetneq E \subseteq K^T$ would be a field $L_i$ for some $i \in \mathbb{N}$ and so $\mathfrak{p}_i \in T$ would not be completely decomposed in $E|K$. Therefore $\hat{T} = D(L|K)$ and $\delta(\hat{T}) = \frac{1}{n} > 0$.

**Remark 4:** There exist infinite sets $T$ of primes such that

$$\delta(\hat{T}) = \delta(T) = 0.$$

*Example 2:* Assume that $K$ is a number field such that
(i) $K$ is not totally real,
(ii) there exists a proper subfield $E$ of $K$ such that $K|E$ is a cyclic extension.

We denote the Galois group $G(K|E)$ by $\Delta$. Let

$$T_0 = \{\mathfrak{p} \,|\, \mathfrak{p} \cap E \text{ is inert in } K|E\}.$$

Then $T_0$ is an infinite set of primes of $K$ of density equal to 0. Let $p$ be a prime number not dividing $[K:E]$ whose extensions to $E$ are completely decomposed or totally ramified in $K|E$. Let $K_\infty$ resp. $E_\infty$ be the compositum of all $\mathbb{Z}_p$-extension of $K$ resp. $E$. Then

$$G(K_\infty|K) = \mathbb{Z}_p^{r_2(K)+1+\delta_K}, \qquad G(E_\infty|E) = \mathbb{Z}_p^{r_2(E)+1+\delta_E},$$

where $r_2$ denotes the number of complex places and $\delta$ is the so-called Leopoldt defect. We have a decomposition of $\mathbb{Z}_p[\Delta]$-modules

$$G(K_\infty|K) \cong G(E_\infty|E) \oplus M,$$

where $M$ is a $\Delta$-module with $M^\Delta = 0$ and $r = \mathrm{rang}_{\mathbb{Z}_p} M \geq r_2(K) - r_2(E) > 0$, since $\delta_K \geq \delta_E$. Let $L$ be the subfield of $K_\infty$ corresponding to $G(K_\infty|K)^\Delta \cong G(E_\infty|E)$, i.e.

$$G(L|K) \cong M$$

and so $G(L|K)^\Delta = 0$. Observe that $L$ is a Galois extension of $E$. Now let $\mathfrak{p} \in T_0$. For the decomposition group with respect to $\mathfrak{p}$ we have the split exact sequence

$$0 \longrightarrow G_\mathfrak{p}(L|K) \longrightarrow G_\mathfrak{p}(L|E) \longrightarrow \Delta_\mathfrak{p} \longrightarrow 0.$$

If $G_\mathfrak{p}(L|K)$ would be non-trivial, then $G_\mathfrak{p}(L|E)$ is not abelian since $\Delta_\mathfrak{p} = \Delta \neq 0$ acts non-trivially on $G_\mathfrak{p}(L|K) \subseteq G(L|K)$. But $\mathfrak{p}$ lies not above $p$ and so it is unramified in $L|K$, thus unramified in $L|E$ and so $G_\mathfrak{p}(L|E)$ is cyclic. Therefore all primes in $T_0$ are completely decomposed in the infinite extension $L|K$. Let

$$T = D(L|K).$$

Then $T$ is an infinite saturated set (containing $T_0$), i.e. $T = \hat{T}$, and since $K^T$ is an infinite extension (containing $L$), it follows from proposition (1.6) that $\delta(\hat{T}) = 0$.

Furthermore, if $S_\mathbb{R} \subseteq T$, then

$$(K^T)_\mathfrak{p} | K_\mathfrak{p} \quad \text{is an infinite extension if } \mathfrak{p} \notin T,$$

hence $T$ is stably saturated. This follows from the fact, that for a prime $\mathfrak{p} \notin D(L|K)$ the non-trivial decomposition group $G_\mathfrak{p}(L|K)$ has to be a infinite subgroup of $G(L|K) \cong \mathbb{Z}_p^r$.

An example for the situation above is a CM-field $K$ with maximal totally real subfield $K^+$. If $T_0(K^+) = \{\mathfrak{p} \,|\, \mathfrak{p} \text{ is inert in } K\}$, then $T_0 = T_0(K)$ is infinite of density equal to zero. Let $K_{ac}$ be the anti-cyclotomic $\mathbb{Z}_p$-extension of $K$, $p \neq 2$, and $T = D(K_{ac}|K)$. Then $T_0 \subseteq T = \hat{T}$ and so $K^T|K$ is an infinite extension and therefore $\delta(\hat{T}) = 0$.

Now we consider the cardinality of an equivalence class $[T]$ of a set of primes $T$ of $K$. By definition, the saturation $\hat{T}$ of $T$ is the unique maximal element of $[T]$ (with respect to the inclusion). But the example 1 of remark 3 shows that there might be infinitely many different (even pairwise disjoint) minimal elements in the class $[T]$.

**Proposition 1.7**

(i) *If $L$ is a finite Galois extension of $K$, then*

$$\#[L] = 1 \quad \text{and} \quad \#[D(L|K)] = \infty.$$

(ii) *If $T$ is a finite set of primes, then*

$$\#[T] = 1 \quad \text{and} \quad \#[K^T] = \infty.$$

**Proof:** Let $L' \in [L]$. Then $L' \subseteq \hat{L} = L$ by theorem (1.5)(ii). Thus $L'|K$ is finite and so
$$L' = K^{D(L'|K)} = K^{D(L|K)} = L.$$

In order to prove the second assertion of (i), let $L$ be a finite Galois extension of $K$ and let $T \subseteq D(L|K)$ be a subset of density equal to zero. Then
$$L = K^{D(L|K)} = K^{D(L|K)\setminus T},$$

hence every subset $D(L|K)\setminus T$ with $\delta(T) = 0$ is contained in $[D(L|K)]$.

By theorem (1.5)(i), every subset of $T$ is (strongly) saturated. Thus the first assertion of (ii) follows. For the second let $S = T \cup \{\mathfrak{p}\}$, where $\mathfrak{p}$ is some prime not in $T$. By [3](10.5.8) we have for the Galois group of the extension $K^T(p)|K^S(p)$ the isomorphism

$$\underset{\mathfrak{p} \in S \setminus T(K^S(p))}{\text{\Large$*$}} G(K_\mathfrak{p}(p)|K_\mathfrak{p}) \xrightarrow{\sim} G(K^T(p)|K^S(p)),$$

where $p$ is some prime number and $E(p)|K$ denotes the maximal $p$-extension inside a Galois extension $E|K$. In particular, the extension $K^T|K^S$ is infinite. Since every Galois extension $L|K$ with $K^S \subsetneq L \subseteq K^T$ is contained in $[K^T]$, the second assertion of (ii) follows. □

# 2 The maximal $p$-extension $k^T(p)$ of $k$

In the following we will consider $p$-extensions and by a saturated set of primes we always mean a $p$-saturated set (see remark 1 of the first section).

Let $S, T$ be sets of primes of a (not necessarily finite) number field $K$ and let $\mathfrak{p}$ be a prime of $K$. Let $K(p)$ be the maximal $p$-extension of $K$. Mostly we will drop the notion $(p)$. So let

$K_{nr}$ is the maximal unramified $p$-extension of $K$,

$K_S$ is the maximal $p$-extension of $K$ which is unramified outside $S$,

$K^T$ is the maximal $p$-extension of $K$ which is completely decomposed at $T$,

$K_S^T$ is the maximal $p$-extension of $K$ which is unramified outside $S$ and

completely decomposed at $T$,

$G_\mathfrak{p}(K) = G_\mathfrak{p}(K(p)|K) \cong G(K_\mathfrak{p}(p)|K_\mathfrak{p})$ is the decomposition group,

$T_\mathfrak{p}(K) = T_\mathfrak{p}(K(p)|K) \cong T(K_\mathfrak{p}(p)|K_\mathfrak{p})$ is the inertia group with respect to $\mathfrak{p}$.

In the following $k$ will always denote a finite number field. If $K|k$ is an infinite extension of number fields and $S$ a set of primes of $k$, then $S(K)$ denotes the profinite space
$$S(K) = \varprojlim_{k'} S(k') \cup \{*_{k'}\}$$
where $k'$ runs through the finite subextensions of $K|k$ and $S(k') \cup \{*_{k'}\}$ is the one-point compactification of the discrete set $S(k')$ of primes of $k'$ lying above $S$. According to [3] (10.5.8) and (10.5.10) we have

**Theorem 2.1** *If $R' \subseteq R \subseteq S \subseteq S'$ are sets of primes of $k$ such that $\delta(S) = 1$ and $R$ is finite, then the canonical homomorphism*
$$\phi_{S',S}^{R',R}: \mathop{\text{\Large$*$}}_{\mathfrak{p} \in S' \backslash S(k_S^R)} T(k_\mathfrak{p}(p)|k_\mathfrak{p}) * \mathop{\text{\Large$*$}}_{\mathfrak{p} \in R \backslash R'(k_S^R)} G(k_\mathfrak{p}(p)|k_\mathfrak{p}) \longrightarrow G(k_{S'}^{R'}|k_S^R)$$
*is an isomorphism. Furthermore we have the following assertions concerning the (strict) cohomological dimension: Assume that $k$ is totally imaginary if $p = 2$, then*
$$cd_p \, G(k_S^R|k) = scd_p \, G(k_S^R|k) = 2.$$

If $L|K$ is a Galois extension we write $H^i(L|K, A)$ for $H^i(G(L|K), A)$, and for a pro-$p$ group $G$ we put $H^i(G) = H^i(G, \mathbb{Z}/p\mathbb{Z})$. For a $p$-extension $K|k$ we will

use the notation

$$\bigoplus_{\mathfrak{p}\in T(K)}' H^i(G_\mathfrak{p}(k(p)|K)) := \varinjlim_{k'} \bigoplus_{\mathfrak{p}\in T(k')} H^i(G_\mathfrak{p}(k(p)|k'))$$

where $k'$ runs through all finite subextensions of $K|k$.

**Proposition 2.2** *Let $p$ be a prime number and let $T$ be a set of primes of a number field $k$ of density $\delta(T) = 0$. Then the following assertions are equivalent:*

(i) *$T$ is stably saturated.*

(ii) *$cd_p\, G_\mathfrak{p}(k(p)|k^T) \leq 1$ for all $\mathfrak{p} \notin T$.*

(iii) *The canonical map*

$$H^2(k(p)|k^T) \xrightarrow{\sim} \bigoplus_{\mathfrak{p}\in T(k^T)}' H^2(G_\mathfrak{p}(k))$$

   *is an isomorphism.*

(iv) *The group $G((k^T)_S|k^T)$ is free for every $S$ with $\delta(S) = 1$, $S \cap T = \varnothing$.*

**Proof:** (i)⇔(ii): By lemma (1.3) the set $T$ is stably saturated if and only if the extensions $(k^T)_\mathfrak{p}|k_\mathfrak{p}$ are infinite for all primes $\mathfrak{p} \notin T$, i.e. if and only if $cd_p\, G_\mathfrak{p}(k(p)|k^T) \leq 1$, see [3] (7.1.8)(i).

In order to prove (ii)⇔(iii), first observe that $k^T|k$ is an infinite extension as $\delta(T) = 0$. Using the Poitou-Tate theorem, see [3] (8.6.10), (10.4.8), and the Hasse principle, loc. cit. (9.1.16), and passing to the limit over all finite extensions $k'$ inside $k^T|k$, we obtain the isomorphism

$$H^2(k(p)|k^T) \xrightarrow{\sim} \bigoplus_{\mathfrak{p}\in T(k^T)}' H^2(G_\mathfrak{p}(k)) \oplus \bigoplus_{\mathfrak{p}\notin T(k^T)}' H^2(G_\mathfrak{p}(k(p)|k^T)),$$

since $\varinjlim_{k'} H^0(G(k(p)|k'), \mu_p)^\vee = 0$. This shows that (iii) is equivalent to

$$H^2(G_\mathfrak{p}(k(p)|k^T)) = 0 \text{ for all } \mathfrak{p} \notin T(k^T),$$

and so we get (ii)⇔(iii).

From the Poitou-Tate exaxt sequences and the Hasse principle for the extensions $k(p)|k$ and $k_S|k$ (using $\delta(S) = 1$) we get the exact sequence

$$0 \longrightarrow H^2(k_S|k) \longrightarrow H^2(k(p)|k) \longrightarrow \bigoplus_{\mathfrak{p}\notin S} H^2(G_\mathfrak{p}(k)) \longrightarrow 0.$$

Passing to the limit from $k$ to $k^T$, we obtain the exact sequence

$$0 \to H^2((k^T)_S|k^T) \to H^2(k(p)|k^T) \to \bigoplus_{\mathfrak{p} \notin S(k^T)}{}' H^2(G_\mathfrak{p}(k^T)) \to 0.$$

Assuming (ii), we have $H^2(G_\mathfrak{p}(k(p)|k^T)) = 0$ for all $\mathfrak{p} \notin T(k^T)$, and so

$$\bigoplus_{\mathfrak{p} \notin S(k^T)}{}' H^2(G_\mathfrak{p}(k(p)|k^T)) = \bigoplus_{\mathfrak{p} \in T(k^T)}{}' H^2(G_\mathfrak{p}(k(p)|k^T)) = \bigoplus_{\mathfrak{p} \in T(k^T)}{}' H^2(G_\mathfrak{p}(k))$$

for every $S$ with $S \cap T = \varnothing$. It follows that $H^2((k^T)_S|k^T) = 0$, i.e. assertion (iv) holds.

Finally, if $H^2((k^T)_{\bar T}|k^T) = H^2(k_{\bar T}|k^T) = 0$ for $\bar T = \mathcal{P}\backslash T$, then, using again the exact sequence above, assertion (iii) holds. $\square$

There is another situation in which the pro-$p$ group $G((k^T)_S|k^T)$ is free.

**Theorem 2.3** *Let $p$ be a prime number and let $T$ and $S$ are sets of primes of a number field $k$, where $T \cup S \neq \mathcal{P}$ is strongly saturated and $T \cap S = \varnothing$. Then the following holds:*

(i) $(k^T)_S = k_{\bar T}$ *and* $\mathrm{cd}_p G((k^T)_S|k^T) \leq 1$.

(ii) *If $\delta(S) = 0$, then $\mathrm{cd}_p G((k^T)_{nr}|k^T) \leq 1$, i.e. the Galois group of the maximal unramified $p$-extension $(k^T)_{nr}|k^T$ is a free pro-$p$-group.*

*In particular, if $T \neq \mathcal{P}$ is a strongly saturated set, then*

$$\mathrm{cd}_p G((k^T)_{nr}|k^T) \leq 1.$$

**Proof:** Since $T \cup S \neq \mathcal{P}$ is stably saturated, we have $\delta(T) \leq \delta(T \cup S) = 0$, hence $\delta(\bar T) = 1$, where $\bar T = \mathcal{P}\backslash T$. It follows from proposition (2.2) that the group $G(k_{\bar T}|k^T)$ is free. Since $T \cup S$ is strongly saturated, the extension $k^{T \cup S}$ realizes the local extensions for all primes in $\overline{T \cup S} = \bar T \backslash S$, and so the extension $k^T$ has this property. Therefore $k_{\bar T}|k^T$ is completely decomposed by $\overline{T \cup S}$, hence $k_{\bar T} = (k^T)_{\bar T} = (k^T)_S$. This proves (i).

Now assume that $\delta(S) = 0$, hence $\delta(\bar T\backslash S) = 1$. Let $K|k$ be a finite extension inside $k^T$. Using the Poitou-Tate theorem and the Hasse principle, we see that the canonical map

$$H^2(K_{\bar T\backslash S}|K) \longrightarrow H^2(k_{\bar T}|K)$$

is injective. Passing to the limit, it follows that

$$H^2((k^T)_{\bar T\backslash S}|k^T) \longrightarrow H^2(k_{\bar T}|k^T)$$

is injective. Since $(k^T)_{\bar{T}\setminus S} = (k^T)_{nr}$, the desired result follows from (i). □

By theorem (1.5)(i) finite sets are strongly saturated. Now we will show that there are also infinite strongly saturated sets.

**Theorem 2.4** *Let $p$ be an odd prime number and let $k$ be a CM-field containing the group $\mu_p$ of all p-th roots of unity, with maximal totally real subfield $k^+$, i.e. $k = k^+(\mu_p)$ is totally imaginary and $[k : k^+] = 2$. Let*

$$T = \{\mathfrak{p} \,|\, \mathfrak{p} \cap k^+ \text{ is inert in } k|k^+\}.$$

*Then the set $T_0 = T \cup S_p$ of primes of $k$ is strongly saturated. Furthermore, the Galois group $G((k^{T_0})_{nr}|k^{T_0})$ is a free pro-p-group.*

**Proof:** Let
$$S = \{\mathfrak{p} \,|\, \mathfrak{p} \cap k^+ \text{ is decomposed in } k|k^+\}$$

and $S_1 = \mathcal{P}\setminus T_0 = S\setminus S_p$. Let $S_0$ be a subset of $S_1$ invariant under the action of $G(k|k^+)$ such that $V = S_1\setminus S_0$ is finite (observe that $T(k_{\mathfrak{p}}(p)|k_{\mathfrak{p}})$ is cyclic for $\mathfrak{p} \in V$ as $V \cap S_p = \varnothing$). Let $K|k$ be a finite extension inside $k^{T_0}$ being Galois over $k^+$ (observe that $k^{T_0}|k^+$ is a Galois extension as $T_0$ is invariant under $G(k|k^+)$). First we show the following

*Claim:* There exists an abelian (not necessarily finite) $p$-extension $L|K$, which is Galois over $k^+$, central over $k$, unramified by $T_0(K)$, completely decomposed at $S_p$ and ramified at each prime of $V(K)$:

$$T(k_{\mathfrak{p}}(p)|K_{\mathfrak{p}})_{G_{\mathfrak{p}}(k(p)|k)} \subseteq G_{\mathfrak{p}}(L|K) \text{ for all } \mathfrak{p} \in V(k).$$

*Proof :* Consider the group extension

$$1 \longrightarrow G(K^{S_p}_{S_1 \cup S_p}|K^{S_p}_{S_0 \cup S_p}) \longrightarrow G(K^{S_p}_{S_1 \cup S_p}|K) \longrightarrow G(K^{S_p}_{S_0 \cup S_p}|K) \longrightarrow 1.$$

Since $\delta(S_0) = 1$, we have $H^2(G(K^{S_p}_{S_0 \cup S_p}|K), \mathbb{Q}_p/\mathbb{Z}_p) = 0$, see (2.1). In the proof of this claim we write $H^i(E|F)$ for $H^i(G(E|F), \mathbb{Q}_p/\mathbb{Z}_p)$. We obtain an exact sequence

$$0 \to H^1(K^{S_p}_{S_0 \cup S_p}|K) \to H^1(K^{S_p}_{S_1 \cup S_p}|K) \to H^1(K^{S_p}_{S_1 \cup S_p}|K^{S_p}_{S_0 \cup S_p})^{G(K^{S_p}_{S_0 \cup S_p}|K)} \to 0,$$

and so an exact sequence

$$0 \to H^1(K^{S_p}_{S_0 \cup S_p}|K)^{G(K|k)} \to H^1(K^{S_p}_{S_1 \cup S_p}|K)^{G(K|k)} \to H^1(K^{S_p}_{S_1 \cup S_p}|K^{S_p}_{S_0 \cup S_p})^{G(K^{S_p}_{S_0 \cup S_p}|k)}$$

$$\to H^1(K|k, H^1(K^{S_p}_{S_0 \cup S_p}|K)) \to H^1(K|k, H^1(K^{S_p}_{S_1 \cup S_p}|K)).$$

Using again that $H^2(G(K_{\tilde{S}}^{S_p}|K), \mathbb{Q}_p/\mathbb{Z}_p) = 0$ where $\tilde{S} = S_0 \cup S_p$ resp. $\tilde{S} = S_1 \cup S_p$, the Hochschild-Serre spectral sequences

$$E_2^{i,j} = H^i(K|k, H^j(K_{\tilde{S}}^{S_p}|K)) \Rightarrow H^{i+j}(K_{\tilde{S}}^{S_p}|k)$$

show that in the commutative diagram

$$\begin{array}{ccc}
H^1(K|k, H^1(K_{S_0 \cup S_p}^{S_p}|K)) & \longrightarrow & H^1(K|k, H^1(K_{S_1 \cup S_p}^{S_p}|K)) \\
\downarrow d_2^{2,1} & & \downarrow d_2^{2,1} \\
H^3(K|k, H^0(K_{S_0 \cup S_p}^{S_p}|K)) & = & H^3(K|k, H^0(K_{S_1 \cup S_p}^{S_p}|K))
\end{array}$$

the differentials $d_2^{2,1}$ are isomorphisms. Thus we obtain an exact sequence

$$0 \longrightarrow \left(G(K_{S_1 \cup S_p}^{S_p}|K_{S_0 \cup S_p}^{S_p})^{ab}\right)_{G(K_{S_0 \cup S_p}^{S_p}|k)} \longrightarrow$$

$$\left(G(K_{S_1 \cup S_p}^{S_p}|K)^{ab}\right)_{G(K|k)} \longrightarrow \left(G(K_{S_0 \cup S_p}^{S_p}|K)^{ab}\right)_{G(K|k)} \longrightarrow 0.$$

Using (2.1), we obtain the isomorphism

$$\prod_{\mathfrak{p} \in V(k)} T(k_{\mathfrak{p}}(p)|K_{\mathfrak{p}})_{G_{\mathfrak{p}}(k(p)|k)} \stackrel{can}{\twoheadrightarrow} G(K_{S_1 \cup S_p}^{S_p}|K_{S_0 \cup S_p}^{S_p})^{ab}_{G(K_{S_0 \cup S_p}^{S_p}|k)} \subseteq G(K_{S_1 \cup S_p}^{S_p}|K)^{ab}_{G(K|k)}.$$

Thus the abelian extension $L|K$ with $G(L|K) = G(K_{S_1 \cup S_p}^{S_p}|K)^{ab}_{G(K|k)}$ has the desired properties and we proved the claim.

As $L|K$ is central over $k$, $G(k|k^+) = <\sigma>$ acts on $G(L|K)$. Thus $G(L|K) = G(L|K)^+ \oplus G(L|K)^-$, where $G(L|K)^{\pm} = G(L|K)^{\sigma \pm 1}$. Let $L^{\pm}$ be defined by $G(L|L^{\pm}) = G(L|K)^{\pm}$, i.e. $G(L^{\pm}|K) = G(L|K)/G(L|K)^{\mp}$.

For $\mathfrak{q} \in T_0 \backslash S_p$ we have the exact sequence

$$1 \longrightarrow G_{\mathfrak{q}}(L^-|k) \longrightarrow G_{\mathfrak{q}}(L^-|k^+) \longrightarrow G_{\mathfrak{q}}(k|k^+) \longrightarrow 1,$$

where $G_{\mathfrak{q}}(L^-|k) = G_{\mathfrak{q}}(L^-|K)$. Suppose that $G_{\mathfrak{q}}(L^-|K) \neq 1$. Since $G(k|k^+) = G_{\mathfrak{q}}(k|k^+)$ acts non-trivially on $G(L^-|K)$, the group $G_{\mathfrak{q}}(L^-|k^+)$ is non-abelian. On the other hand the extension $L^-|k^+$ is unramified at $\mathfrak{q}$. This contradiction shows that all primes of $T_0 \backslash S_p$ are completely decomposed in $L^-|K$, and so in $L^-|k$. Since $L|K$ is completely decomposed at $S_p$, we obtain $L^- \subseteq k^{T_0}$.

Let $\mathfrak{p} \in V$ and so $\mathfrak{p} \cap k^+$ splits in $k|k^+$. Let $\bar{\mathfrak{p}} = \mathfrak{p}^\sigma$ be the conjugated prime and let $\mathfrak{P}$ be a prolongation of $\mathfrak{p}$ to $K$ and $\bar{\mathfrak{P}} = \mathfrak{P}^\sigma$. By the claim it follows that

$$\left(T(k_{\mathfrak{P}}(p))|K_{\mathfrak{P}})_{G_{\mathfrak{P}}(k(p))|k)} \oplus T(k_{\bar{\mathfrak{P}}}(p))|K_{\bar{\mathfrak{P}}})_{G_{\bar{\mathfrak{P}}}(k(p))|k)}\right)^-$$

injects into $G(L|K)^- \stackrel{\sim}{\to} G(L^-|K) = G(L|K)/G(L|K)^+$, and so $L^-|K$ is ramified by $\mathfrak{P}$.

Varying the set $S_0$ and the extension $K|k$, it follows that the extension $k^{T_0}|k$ realizes the ramified part of $k_{\mathfrak{p}}(p)|k_{\mathfrak{p}}$ for all $\mathfrak{p} \in S_1 = \mathcal{P} \backslash T_0$. Since $k^{T_0}|k$ is a Galois extension, also the unramified part must be realized, i.e.

$$(k^{T_0})_{\mathfrak{P}} = k_{\mathfrak{p}}(p) \quad \text{for all } \mathfrak{P}|\mathfrak{p},\ \mathfrak{p} \notin T_0$$

(if $k_{\mathfrak{P}}(p)|(k^{T_0})_{\mathfrak{P}}$ would have a non-trivial unramified part, then, as the subgroup generated by the Frobenius automorphism is not normal, this extension would also have a ramified part).

The last assertion of the theorem follows from theorem (2.3). □

Let

$$\coprod_{t \in T}(A_t, B_t) = \{(a_t)_{t \in T} \in \prod_{t \in T} A_t \mid a_t \in B_t \text{ for almost all } t \in T\}$$

be the restricted product over a discrete set $T$ of abelian locally compact groups $A_t$ with respect to closed subgroups $B_t$. The topology is given by the subgroups $V$ with

(i) $V \cap A_t$ is open in $A_t$ for all $t \in T$,
(ii) $V \supseteq B_t$ for almost all $t \in T$.

Then we call

$$\coprod_{t \in T}^{c}(A_t, B_t) := \varprojlim_{V} \left(\coprod_{t \in T}(A_t, B_t)\right)/V,$$

the *compactification* of $\coprod_{t \in T}(A_t, B_t)$, where $V$ runs through all open subgroups of finite index in $\coprod_{t \in T}(A_t, B_t)$. The the canonical map $\coprod_{t \in T}(A_t, B_t) \to \coprod_{t \in T}^{c}(A_t, B_t)$ has dense image.

We define the *discretization* of $\coprod_{t \in T}(A_t, B_t)$ by

$$\coprod_{t \in T}^{d}(A_t, B_t) := \varinjlim_{W} W$$

where $W$ runs through the finite subgroups of $\coprod_{t \in T}(A_t, B_t)$. If the subgroups $B_t$ of $A_t$, $t \in T$, are open and compact, then $\coprod_{t \in T}(A_t, B_t)$ is locally compact. Using the equality

$$(\coprod_{t \in T}(A_t, B_t))^{\vee} = \coprod_{t \in T}(A_t^{\vee}, (A_t/B_t)^{\vee}),$$

we obtain

$$\coprod_{t \in T}^{d}(A_t, B_t) = (((\coprod_{t \in T}(A_t, B_t))^{\vee})^c)^{\vee} \quad \text{and} \quad \coprod_{t \in T}^{c}(A_t, B_t) = (((\coprod_{t \in T}(A_t, B_t))^{\vee})^d)^{\vee},$$

where $^{\vee}$ denotes the Pontryagin-dual.

**Proposition 2.5** *Let $T$ be a set of primes of a number field $k$ with $\delta(T) = 0$. Then*

(i) $G(k(p)/k^T)/G(k(p)/k^T)^* = \varprojlim_{K} \prod^c_{\mathfrak{P} \in T(K)} (G_{\mathfrak{P}}(k)/G_{\mathfrak{P}}(k)^*, \tilde{T}_{\mathfrak{P}}(k))$,

(ii) $H^1(G(k(p)/k^T)) = \varinjlim_{K} \prod^d_{\mathfrak{P} \in T(K)} (H^1(G_{\mathfrak{P}}(k)), H^1_{nr}(G_{\mathfrak{P}}(k)))$,

*where $K$ runs through the finite subextensions of $k^T|k$ and $\tilde{T}_{\mathfrak{P}}(k)$ is the group $T_{\mathfrak{P}}(k)G_{\mathfrak{P}}(k)^*/G_{\mathfrak{P}}(k)^*$.*

**Proof:** Since the set $\mathcal{P}\backslash T$ has density equal to 1, we get from the Poitou-Tate exact sequence and the Hasse principle the commutative and exact diagram

$$\begin{array}{ccccccc}
0 & \longrightarrow & H^1(\bar{k}|k, \mu_p) & \longrightarrow & \prod_{\mathfrak{p} \notin T} H^1(\bar{k}_{\mathfrak{p}}|k_{\mathfrak{p}}, \mu_p) & & \\
& & \| & & \uparrow & & \\
0 & \longrightarrow & H^1(\bar{k}|k, \mu_p) & \longrightarrow & \prod_{\mathfrak{p}} H^1(\bar{k}_{\mathfrak{p}}|k_{\mathfrak{p}}, \mu_p) & \longrightarrow H^1(\bar{k}|k, \mathbb{Z}/p\mathbb{Z})^{\vee} \longrightarrow 0 \\
& & & & \uparrow & & \uparrow \\
& & & & \prod_{\mathfrak{p} \in T} H^1(\bar{k}_{\mathfrak{p}}|k_{\mathfrak{p}}, \mu_p) & \xrightarrow{\sim} \prod_{\mathfrak{p} \in T} H^1(\bar{k}_{\mathfrak{p}}|k_{\mathfrak{p}}, \mathbb{Z}/p\mathbb{Z})^{\vee}
\end{array}$$

where we use the local duality theorem

$$H^1(\bar{k}_{\mathfrak{p}}|k_{\mathfrak{p}}, \mu_p) \xrightarrow{\sim} H^1(G_{\mathfrak{p}}(k), \mathbb{Z}/p\mathbb{Z})^{\vee} = G_{\mathfrak{p}}(k)/G_{\mathfrak{p}}(k)^*$$

and for $\mathfrak{p}$ not above $p$ or $\infty$ the isomorphism

$$H^1_{nr}(\bar{k}_{\mathfrak{p}}|k_{\mathfrak{p}}, \mu_p) \xrightarrow{\sim} (H^1(T_{\mathfrak{p}}(k), \mathbb{Z}/p\mathbb{Z})^{G_{\mathfrak{p}}(k)})^{\vee} \cong T_{\mathfrak{p}}(k)G_{\mathfrak{p}}(k)^*/G_{\mathfrak{p}}(k)^*,$$

see [3] (8.6.10), (9.1.9), (9.1.10), (7.2.15). Thus we get a continuous injection

$$\prod_{\mathfrak{p} \in T(k)} (G_{\mathfrak{p}}(k)/G_{\mathfrak{p}}(k)^*, \tilde{T}_{\mathfrak{p}}(k)) \hookrightarrow G(k(p)/k)/G(k(p)/k)^*$$

and so an injection $\prod^c_{\mathfrak{p} \in T(k)}(G_{\mathfrak{p}}(k)//G_{\mathfrak{p}}(k)^*, \tilde{T}_{\mathfrak{p}}(k)) \hookrightarrow G(k(p)/k)/G(k(p)/k)^*$. Passing to the limit, we obtain the injection

$$\varprojlim_{K} \prod^c_{\mathfrak{P} \in T(K)} (G_{\mathfrak{P}}(k)/G_{\mathfrak{P}}(k)^*, \tilde{T}_{\mathfrak{P}}(k)) \hookrightarrow G(k(p)/k^T)/G(k(p)/k^T)^*$$

which, by definition of the field $k^T$, is also surjective. Therefore we obtain (i) and (ii) is just the dual assertion. □

**Proposition 2.6** *Let $p$ be a prime number and let $T \neq \mathcal{P}$ be a stably saturated set of primes of a number field $k$. Assume that $k$ is totally imaginary if $p = 2$. Then*
$$cd_p \, G(k^T|k) \leq 2.$$

**Proof:** We consider the Hochschild-Serre spectral sequence
$$E_2^{ij} = H^i(G/N, H^j(N)) \Rightarrow E^{i+j} = H^{i+j}(G)$$
for an extension $1 \to N \to G \to G/N \to 1$ of pro-$p$ groups. If $cd_p(G) \leq 2$, $E_2^{11} = 0$ and the edge homomorphism $E^2 \to E_2^{02}$ is surjective, then it follows that $E_2^{30} = 0$, i.e. $cd_p(G/N) \leq 2$.

Now let $G = G(k(p)|k)$ and $N = G(k(p)|k^T)$. Using (2.2) and the Poitou-Tate theorem, it follows that the canonical map
$$H^2(k(p)|k) \twoheadrightarrow \bigoplus_{\mathfrak{p} \in T(k)} H^2(G_\mathfrak{p}(k)) \xrightarrow{\sim} H^0(k^T|k, H^2(k(p)|k^T))$$
is surjective. Furthermore, using (2.5)(ii), we see that

$$\begin{aligned}
H^1(k(p)|k^T) &= \varinjlim_K H^1(\coprod^c_{\mathfrak{p} \in T(K)}(G_\mathfrak{p}(k)/G_\mathfrak{p}(k)^*, \tilde{T}_\mathfrak{p}(k))) \\
&= \varinjlim_K Coind^{G(K|k)} H^1(\coprod^c_{\mathfrak{p} \in T(k)}(G_\mathfrak{p}(k)/G_\mathfrak{p}(k)^*, \tilde{T}_\mathfrak{p}(k))) \\
&= Coind^{G(k^T|k)} H^1(\coprod^c_{\mathfrak{p} \in T(k)}(G_\mathfrak{p}(k)/G_\mathfrak{p}(k)^*, \tilde{T}_\mathfrak{p}(k)))
\end{aligned}$$

is a $G(k^T|k)$-coinduced module, hence $H^1(k^T|k, H^1(k(p)|k^T)) = 0$. Therefore we obtain the desired result. $\square$

If $T$ is a set of primes of $k$, then let $T(k^T) = \varprojlim_K T(K) \cup \{*_K\}$, where $K$ runs through the finite subextensions of $k^T/k$. The following theorem asserts that the Galois group $G(k(p)/k^T)$, where $T \neq \mathcal{P}$ is stably saturated, is a corestricted free pro-$p$-product of the family $(G_\mathfrak{P}(k))_{\mathfrak{P} \in T(k^T)}$ of decomposition groups with respect to the continuous family $(T_\mathfrak{P}(k))_{\mathfrak{P} \in T(k^T)}$ of inertia groups, see [2] for the definition.

**Theorem 2.7** *Let $T \neq \mathcal{P}$ be a stably saturated set of primes of a number field $k$. Then the canonical map*
$$\varphi : \mathop{\text{\Large $*$}}_{\mathfrak{P} \in T(k^T)}(G_\mathfrak{P}(k), T_\mathfrak{P}(k)) \xrightarrow{\sim} G(k(p)/k^T)$$
*is an isomorphism.*

**Proof:** We can apply [2] prop. (4.3) and have to show that the induced maps

$$\varphi_* : H^1(G(k(p)/k^T), \mathbb{Z}/p\mathbb{Z}) \xrightarrow{\sim} \varinjlim_K \prod_{T(K)}^d (H^1(G_{\mathfrak{P}}(k), \mathbb{Z}/p\mathbb{Z}), H^1_{nr}(G_{\mathfrak{P}}(k), \mathbb{Z}/p\mathbb{Z}))$$

$$\varphi_* : H^2(G(k(p)/k^T), \mathbb{Z}/p\mathbb{Z}) \hookrightarrow \varinjlim_K \prod_{T(K)}^d (H^2(G_{\mathfrak{P}}(k), \mathbb{Z}/p\mathbb{Z}), H^2_{nr}(G_{\mathfrak{P}}(k), \mathbb{Z}/p\mathbb{Z}))$$

are bijective resp. injective. Since $H^2_{nr}(G_{\mathfrak{P}}(k)) = 0$, it follows that

$$H^2(G) = \varinjlim_K \prod_{T(K)}^d (H^2(G_{\mathfrak{P}}(k)), H^2_{nr}(G_{\mathfrak{P}}(k))) = \bigoplus_{\mathfrak{p} \in T(k^T)}' H^2(G_{\mathfrak{p}}(k)).$$

Now the result follows from (2.5)(ii) and (2.2). □

In the situation of theorem (2.4) we obtain

**Corollary 2.8** *Let $p$ be an odd prime number and let $k$ be a CM-field containing the group $\mu_p$ of all $p$-th roots of unity, with maximal totally real subfield $k^+$, i.e. $k = k^+(\mu_p)$ is totally imaginary and $[k : k^+] = 2$. Let*

$$T = \{\mathfrak{p} \mid \mathfrak{p} \cap k^+ \text{ is inert in } k|k^+\} \cup S_p.$$

*Then*

$$\underset{\mathfrak{P} \in T(k^T)}{\text{\Huge{*}}} (G_{\mathfrak{P}}(k), T_{\mathfrak{P}}(k)) \xrightarrow{\sim} G(k(p)/k^T).$$

Mathematisches Institut
der Universität Heidelberg
Im Neuenheimer Feld 288
69120 Heidelberg
Germany
e-mail: wingberg (at) mathi.uni-heidelberg.de